\input amstex
\documentstyle{amsppt}
\magnification=\magstep1                        %<====
\hsize6.5truein\vsize8.9truein                  %<====
\NoRunningHeads
\loadeusm

\magnification=\magstep1                        %<====
\hsize6.5truein\vsize8.9truein                  %<====
\NoRunningHeads
\loadeusm

\document
\topmatter

\title
The number of unimodular zeros of self-reciprocal polynomials with coefficients 
in a finite set
\endtitle

\rightheadtext{the multiplicity of the zero at $1$ of polynomials}

\author Tam\'as Erd\'elyi
\endauthor

%\address  Department of Mathematics and Statistics,
%Simon Fraser University, Burnaby, B.C., Canada V5A 1S6 (P. Borwein) \endaddress

\address Department of Mathematics, Texas A\&M University,
College Station, Texas 77843, College Station, Texas 77843 (T. Erd\'elyi) \endaddress

%\address Mathematical Institute, L\'or\'and E\"otv\"os University, P\'azm\'any P. s. 1/c, Budapest, Hungary H-1117; 
%and Computer and Automation Research Institute, Kende u. 13-17. Budapest, Hungary H-1111
%(G. K\'os) \endaddress

\thanks {{\it 2010 Mathematics Subject Classifications.} 11C08, 41A17, 26C10, 30C15}
\endthanks

\keywords
self-reciprocal polynomials, trigonometric polynomials, restricted coefficients, number of zeros on the unit circle,
number of real zeros in a period, Conrey's question
\endkeywords

\date February 4, 2016
\enddate

\abstract
We study the number $\text {\rm NZ}(T_n)$ of real zeros of trigonometric polynomials
$$T_n(t) = \sum_{j=0}^n{a_{j,n} \cos(jt)}$$
in a period $[a,a+2\pi)$, $a \in {\Bbb R}$, and the number $\text {\rm NZ}(P_n)$ of zeros of 
self-reciprocal algebraic polynomials
$$P_n(z) = \sum_{j=0}^n{a_{j,n} z^j}$$
on the unit circle. One of the highlights of this paper states that 
$\lim_{n \rightarrow \infty}{\text {\rm NZ}(T_n)} = \infty$ whenever the set 
$$\{a_{j,n}: j=0,1,\ldots,n, \enskip n \in {\Bbb N}\} \subset [0,\infty)$$
is finite and 
$$\lim_{n \rightarrow \infty}{|\{j \in \{0,1,\ldots,n\}: a_{j,n} \neq 0\}|} = \infty\,.$$
This follows from a more general result stating that $\lim_{n \rightarrow \infty}{\text {\rm NZ}(T_n)} = \infty$
whenever the set 
$$\{a_{j,n}: j=0,1,\ldots,n, \enskip n \in {\Bbb N}\} \subset {\Bbb R}$$
is finite and $\lim_{n \rightarrow \infty}{\text {\rm NC}_k(T_n)} = \infty$ for every $k \in {\Bbb N}$, where 
$$\text {\rm NC}_k(T_n) := \left|\left\{u: \enskip 0 \leq u \leq n-k+1,\sum_{j=u}^{u+k-1}{a_{j,n}} \neq 0 \right\}\right|\,.$$
\endabstract
\endtopmatter

\document

\head 1. Introduction and Notation \endhead
Research on the distribution of the zeros of algebraic polynomials has a 
long and rich history. In fact all the papers [1-43] in our list of references 
are just some of the papers devoted to this topic. The study of the number of 
real zeros trigonometric polynomials and the number of unimodular zeros 
(that is, zeros lying on the unit circle of the complex plane) of algebraic 
polynomials with various constraints on their coefficients are the subject of 
quite a few of these. We do not try to survey these in our introduction.

Let ${\Cal P}_n^c$ denote the set of all algebraic polynomials of degree at most $n$ with complex coefficients.
Let $S \subset {\Bbb C}$. Let ${\Cal P}_n^c(S)$ be the set of all algebraic polynomials of
degree at most $n$ with each of their coefficients in $S$. A polynomial
$$P_n(z) = \sum_{j=0}^n{a_jz^j} \tag 1.1$$
is called conjugate-reciprocal if
$$\overline{a}_j = a_{n-j}, \qquad j=0,1,\ldots,n\,.$$
A polynomial $P_n$ of the form (1.1) is called plain-reciprocal or self-reciprocal if
$$a_j = a_{n-j}, \qquad j=0,1,\ldots,n\,.$$
If a conjugate reciprocal polynomial $P_n$ has only real coefficients, then it is obviously
plain-reciprocal. Associated with an algebraic polynomial
$$P_n(z) = \sum_{j=0}^n{a_{j,n}z^j}$$
we introduce the numbers
$$\text{\rm NC}(P_n) := |\{j \in \{0,1,\ldots,n\}: a_{j,n} \neq 0\}|\,.$$
Let $\text{\rm NZ}(P_n)$ denote the number of real zeros (by counting multiplicities) of an algebraic
polynomial $P_n$ on the unit circle. Associated with a trigonometric polynomial
$$T_n(t) = \sum_{j=0}^n{a_{j,n} \cos(jt)}$$
we introduce the numbers
$$\text{\rm NC}(T_n) := |\{j \in \{0,1,\ldots,n\}: a_{j,n} \neq 0\}|\,.$$
Let $\text{\rm NZ}(T_n)$ denote the number of real zeros (by counting multiplicities) of a trigonometric
polynomial $T_n$ in a period $[a,a+2\pi)$, $a \in {\Bbb R}$. The quotation below is from [6]. 

``Let $0 \leq n_1 < n_2 < \cdots < n_N$ be integers. A cosine polynomial of the form
$T_N(\theta) = \sum_{j=1}^N {\cos (n_j\theta)}$ must have at least one real zero in
a period $[a,a+2\pi)$, $a \in {\Bbb R}$. This is obvious if $n_1 \neq 0$, since 
then the integral of the sum on a period is $0$.
The above statement is less obvious if $n_1 = 0$, but for
sufficiently large $N$ it follows from Littlewood's Conjecture simply. Here we mean the
Littlewood's Conjecture proved by S. Konyagin [25] and independently by McGehee, Pigno,
and Smith [33] in 1981. See also [13, pages 285-288] for a book proof.
It is not difficult to prove the statement in general even in the case $n_1=0$ without 
using Littlewood's Conjecture. One possible way is to use the identity
$$\sum_{j=1}^{n_N}{T_N((2j-1)\pi/n_N)} = 0\,.$$
See [26], for example.
Another way is to use Theorem 2 of [34]. So there is certainly no shortage of possible
approaches to prove the starting observation of this paper even in the case $n_1=0$.

It seems likely that the number of zeros of the above sums in a period
must tend to $\infty$ with $N$. In a private communication B. Conrey asked
how fast the number of real zeros of the above sums in a period tends to $\infty$ as
a function $N$. In [4] the authors observed that
for an odd prime $p$ the Fekete polynomial $f_p(z)=\sum^{p-1}_{k=0} \big({k\over p}\big) z^k$
(the coefficients are Legendre symbols) has $\sim \kappa_0 p$  zeros on the unit circle, where
$0.500813>\kappa_0>0.500668$. Conrey's question in general does not appear to be easy.

Littlewood  in his 1968 monograph ``Some Problems in Real and Complex Analysis"
[10, problem~22] poses the following research problem,
which appears to still be open: ``If the $n_m$ are integral and all
different, what is the lower bound on the number of real zeros of
$\sum_{m=1}^N \cos (n_m\theta)$?  Possibly $N-1$, or not much less." 
Here real zeros are counted in a period.
In fact no progress appears to have been made on this
in the last half century. In a recent paper [3] we showed that this is
false. There exists a cosine polynomials $\sum_{m=1}^N \cos (n_m\theta)$ with
the $n_m$ integral and all different so that the number of its real zeros
in the period is $O(N^{9/10}(\log N)^{1/5})$ (here the frequencies $n_m = n_m(N)$
may vary with $N$). However, there
are reasons to believe that a cosine polynomial $\sum_{m=1}^N \cos (n_m\theta)$
always has many zeros in the period." 

One of the highlights of this paper is to show that the number of real zeros of the sums
$T_N(\theta) = \sum_{j=1}^N {\cos (n_j\theta)}$ in a period tends to $\infty$ whenever 
$0 \leq n_1 < n_2 < \cdots < n_N$ are integers and $N$ tends to $\infty$, even though  
the part "how fast" in Conrey's question remains open. In fact, we will prove more general 
results of this variety. Let 
$${\Cal L}_n := \left\{P: \enskip P(z)=\sum_{j=0}^n{a_jz^j}: \enskip a_j \in \{-1,1\} \right\} \,.$$ 
Elements of ${\Cal L}_n$ are often called Littlewood polynomials of degree $n$. Let 
$${\Cal K}_n := \left\{P: \enskip P(z)=\sum_{j=0}^n{a_jz^j}: \enskip a_j \in {\Bbb C}, \enskip |a_0|=|a_n|=1,
\enskip |a_j| \leq 1 \right\}\,,$$
Observe that ${\Cal L}_n \subset {\Cal K}_n$. 
In [10] we proved that any polynomial $P \in {\Cal K}_n$
has at least $8n^{1/2}\log n$ zeros in any open disk centered at a point on the unit circle 
with radius $33n^{-1/2}\log n$. Thus polynomials in ${\Cal K}_n$
have a few zeros near the unit circle. One may naturally ask how many unimodular roots a 
polynomial in ${\Cal K}_n$ can have. 
Mercer [34] proved that if a Littlewood polynomial $P \in {\Cal L}_n$ of the form (1.1) 
is skew reciprocal, that is, $a_j = (-1)^ja_{n-j}$ for each $j=0,1,\ldots,n$, then it has 
no zeros on the unit circle. However, by using different elementary methods it was observed 
in both [18] and [34] that if a Littlewood polynomial $P$ of the form (1.1) is self-reciprocal, 
that is  $a_j = a_{n-j}$ for each $j=0,1,\ldots,n$, $n \geq 1$, then it has at least one zero 
on the unit circle. 
Mukunda [35] improved this result by showing that every self-reciprocal Littlewood polynomial of 
odd degree at least $3$ has at least $3$ zeros on the unit circle. Drungilas [16] proved that 
every self-reciprocal Littlewood polynomial of odd degree $n \geq 7$ has at least $5$ zeros 
on the unit circle and every self-reciprocal Littlewood polynomial of even degree $n \geq 14$ 
has at least $4$ zeros on the unit circle. In [4] two types of Littlewood polynomials are considered: 
Littlewood polynomials with one sign change in the sequence of coefficients and Littlewood polynomials 
with one negative coefficient, and the numbers of the zeros such Littlewood polynomials have   
on the unit circle and inside the unit disk, respectively, are investigated. Note that the 
Littlewood polynomials studied in [4] are very special. In [7] we proved that the average number 
of zeros of self-reciprocal Littlewood polynomials of degree $n$ is at least $n/4$. However, 
it is much harder to give decent lower bounds for the quantities 
$$\text {\rm NZ}_n := \min_{P}{\text {\rm NZ}(P)}\,,$$
where $\text {\rm NZ}(P)$ denotes the number of zeros of a polynomial $P$ lying on the
unit circle and the minimum is taken for all self-reciprocal Littlewood polynomials 
$P \in {\Cal L}_n$. It has been conjectured for a long time that 
$\lim_{n \rightarrow \infty}\text {\rm NZ}_n = \infty$.  
In this paper we show that 
$\lim_{n \rightarrow \infty}\text {\rm NZ}(P_n) = \infty$  
whenever $P_n \in {\Cal L}_n$ is self-reciprocal and 
$\lim_{n \rightarrow \infty}{|P_n(1)|} = \infty$.
This follows as a consequence of a more general result in which the 
coefficients of the self-reciprocal polynomials $P_n$ of degree at most $n$ 
belong to a fixed finite set of real numbers. In [6] we proved the following result.

\proclaim{Theorem 1.1} If the set $\{a_j:j \in {\Bbb N}\} \subset {\Bbb R}$ is finite, the set
$\{j \in {\Bbb N}: a_j \neq 0\}$ is infinite, the sequence  $(a_j)$ is not eventually periodic, and 
$$T_n(t) = \sum_{j=0}^n{a_j \cos(jt)}\,,$$
then $\lim_{n \rightarrow \infty}{\text {\rm NZ}(T_n)} = \infty\,.$ 
\endproclaim

In [6] Theorem 1.1 is stated without the assumption that the sequence 
$(a_j)$ is not eventually periodic. However, as the following example shows, 
Lemma 3.4 in [6], dealing with the case of eventually periodic sequences $(a_j)$, 
is incorrect. Let 
$$\split T_n(t) := & \cos t + \cos((4n+1)t) + \sum_{k=0}^{n-1}{(\cos((4k+1)t) - \cos((4k+3)t))} \cr 
= & \frac{1+\cos((4n+2)t)}{2\cos t} + \cos t\,. \cr \endsplit$$
It is easy to see that $T_n(t) \neq 0$ on $[-\pi,\pi] \setminus \{-\pi/2,\pi/2\}$ and the 
zeros of $T_n$ at $-\pi/2$ and $\pi/2$ are simple. Hence $T_n$ has only two (simple) zeros 
in the period. So the conclusion of Theorem 1.1 above is false for  
the sequence $(a_j)$ with $a_0 := 0$, $a_1 := 2$, $a_3 :=-1$, $a_{2k} := 0$, $a_{4k+1} := 1$, $a_{4k+3} := -1$ 
for every $k=1,2,\ldots$. Nevertheless, Theorem 1.1 can be saved even 
in the case of eventually periodic sequences $(a_j)$ if we 
assume that $a_j \neq 0$ for all sufficiently large $j$. See Lemma 3.11.  
So Theorem 1 in [6] can be corrected as 

\proclaim{Theorem 1.2} If the set $\{a_j:j \in {\Bbb N}\} \subset {\Bbb R}$ is finite, 
$a_j \neq 0$ for all sufficiently large $j$, and   
$$T_n(t) = \sum_{j=0}^n{a_j \cos(jt)}\,,$$
then $\lim_{n \rightarrow \infty}{\text {\rm NZ}(T_n)} = \infty\,.$
\endproclaim

It was expected that the conclusion of the above theorem remains true 
even if the coefficients of $T_n$ do not come from the same sequence, 
that is 
$$T_n(t) = \sum_{j=0}^n{a_{j,n} \cos(jt)}\,,$$
where the set 
$$\{a_{j,n}:j=0,1,\ldots,n, n \in {\Bbb N}\} \subset {\Bbb R}$$ 
is finite and 
$$\lim_{n \rightarrow \infty}{|\{j: a_{j,n} \neq 0\}|} = \infty\,.$$
The purpose of this paper is to prove such an extension of Theorem 1.1. 
The already mentioned Littlewood Conjecture, proved by Konyagin [25] and
independently by McGehee, Pigno, and B. Smith [33], plays an key role 
in the proof of the main results in this paper. This states the following.

\proclaim{Theorem 1.3}
There is an absolute constant $c > 0$ such that
$$\int_{0}^{2\pi}{\Big| \sum_{j=1}^m{a_je^{i\lambda_jt}} \Big| \, dt} \geq c\gamma \log m\,$$
whenever $\lambda_1, \lambda_2, \ldots, \lambda_n$ are distinct integers and 
$a_1, a_2, \ldots, a_m$ are complex numbers of modulus at least $\gamma > 0$. 
\endproclaim

This is an obvious consequence of the following result a book proof of which has been 
worked out by Lorentz in [13, pages 285-288].

\proclaim{Theorem 1.4}
If $\lambda_1 < \lambda_2 < \cdots < \lambda_m$ are integers and $a_1, a_2, \ldots, a_m$ are complex numbers, then 
$$\int_{0}^{2\pi}{\Big| \sum_{j=1}^m{a_je^{i\lambda_j t}} \Big| \, dt} \geq
\frac{1}{30} \sum_{j=1}^m{\frac{|a_j|}{j}}\,.$$
\endproclaim

\head 2. New Results \endhead

Associated with an algebraic polynomial 
$$P_n(t) = \sum_{j=0}^n{a_{j,n} z^j}\,, \qquad a_{j,n} \in {\Bbb C}\,,$$
let
$$\text {\rm NC}_k(P_n) := \left|\left\{u: \enskip 0 \leq u \leq n-k+1,\sum_{j=u}^{u+k-1}{a_{j,n}} \neq 0 \right\}\right|\,.$$

\proclaim{Theorem 2.1}
If $S \subset {\Bbb R}$ is a finite set, $P_{2n} \in {\Cal P}_{2n}^c(S)$ are self-reciprocal polynomials,
$$T_n(t) := P_{2n}(e^{it})e^{-int}\,,$$ 
and
$$\lim_{n \rightarrow \infty}{\text {\rm NC}_k(P_{2n})} = \infty \tag 2.1$$
for every $k \in {\Bbb N}$, then
$$\lim_{n \rightarrow \infty}{\text {\rm NZ}(T_n)} = \infty\,. \tag 2.2$$
\endproclaim

\proclaim{Corollary 2.2}
If $S \subset {\Bbb R}$ is a finite set, $P_{2n} \in {\Cal P}_{2n}^c(S)$ are self-reciprocal polynomials, 
$$T_n(t) := P_{2n}(e^{it})e^{-int}\,,$$ 
and
$$\lim_{n \rightarrow \infty}{|P_{2n}(1)|} = \infty\,, \tag 2.3$$
then (2.2) holds.
\endproclaim

Our next result is slightly more general than Corollary 2.2, and it follows from Corollary 2.2 simply.  

\proclaim{Corollary 2.3}
If $S \subset {\Bbb R}$ is a finite set, $P_n \in {\Cal P}_n^c(S)$ are self-reciprocal polynomials,
and
$$\lim_{n \rightarrow \infty}{|P_n(1)|} = \infty\,, \tag 2.4$$
then
$$\lim_{n \rightarrow \infty}{\text {\rm NZ}(P_n)} = \infty\,. \tag 2.5$$
\endproclaim

We say that $S \subset {\Bbb R}$ has property (2.6) if (for every $k \in {\Bbb N}$)  
$$s_1 + s_2 + \cdots + s_k = 0\,, \enskip s_1,s_2,\ldots,s_k \in S\,, 
\enskip \text {\rm implies} \enskip s_1=s_2= \cdots =s_k=0\,,\tag 2.6$$ 
that is, any sum of nonzero elements of $S$ is different from zero.   

\proclaim{Corollary 2.4}
If the finite set $S \subset {\Bbb R}$ has property (2.6), $P_{2n} \in {\Cal P}_{2n}^c(S)$ are self-reciprocal polynomials,  
$$T_n(t) := P_{2n}(e^{it})e^{-int}\,,$$ 
and 
$$\lim_{n \rightarrow \infty}{\text {\rm NC}(P_{2n})} = \infty\,, \tag 2.7$$
then (2.2) holds. 
\endproclaim

Our next result is slightly more general than Corollary 2.4, and it follows from Corollary 2.4 simply.

\proclaim{Corollary 2.5}
If the finite set $S \subset {\Bbb R}$ has property (2.5),
$P_n \in {\Cal P}_n^c(S)$ are self-reciprocal polynomials, and 
$$\lim_{n \rightarrow \infty}{\text {\rm NC}(P_n)} = \infty\,, \tag 2.8$$
then (2.5) holds. 
\endproclaim

Our next result is an obvious consequence of Corollary 2.2. 

\proclaim{Corollary 2.6}
If 
$$T_n(t) = \sum_{j=0}^n{a_{j,n} \cos(jt)}\,,$$
where the set
$$S := \{a_{j,n}: j=0,1,\ldots,n, \enskip n \in {\Bbb N}\} \subset {\Bbb R}$$
is finite and 
$$\lim_{n \rightarrow \infty}{\Big| \sum_{j=0}^n{a_{j,n}} \Big|} = \infty\,,$$
then (2.2) holds.
\endproclaim

Our next result is an obvious consequence of Corollary 2.6.

\proclaim{Corollary 2.7}
If 
$$T_n(t) = \sum_{j=0}^{\infty}{a_{j,n} \cos(jt)}\,,$$
where the set
$$S := \{a_{j,n}: j=0,1,\ldots,n, \enskip n \in {\Bbb N}\} \subset [0,\infty)$$
is finite, and  
$$\lim_{n \rightarrow \infty}{\text{\rm NC}(T_n)} = \infty\,,$$ 
then (2.2) holds.
\endproclaim

\head 3. Lemmas \endhead

In our first seven Lemmas we assume that $S \subset {\Bbb C}$ is a finite set, $P_{2n} \in {\Cal P}_{2n}^c(S)$,  
and we use the notation 
$$T_n(t) := P_{2n}(e^{it})e^{-int}\,.$$

\proclaim {Lemma 3.1}
If $S \subset {\Bbb C}$ is a finite set, $P_{2n} \in {\Cal P}_{2n}^c(S)$,  and $H \in {\Cal P}_m^c$ is a polynomial of minimal degree $m$ 
such that
$$\sup_{n \in {\Bbb N}}{\text {\rm NC}(P_{2n}H)} < \infty\,, \tag 3.1$$
then each zero of $H$ is a root of unity, and each zero of $H$ is simple.
\endproclaim

\demo{Proof}
Let $H \in {\Cal P}_m^c$ satisfy the assumptions of the lemma and suppose to the contrary that 
$H(\alpha) = 0$, where $0 \neq \alpha \in {\Bbb C}$ is not a root of unity.
Let $G \in {\Cal P}_{m-1}^c$ be defined by
$$G(z) := \frac{H(z)}{z-\alpha}\,. \tag 3.2$$
Let $S_n^*$ be the set of the coefficients of $P_{2n}G$, and let 
$$S^* := \bigcup_{n \in {\Bbb N}}{S_n^*}\,.$$ 
As $P_{2n} \in {\Cal P}_{2n}(S)$ and the set $S$ is finite, the set $S^*$ is also finite.
Let
$$(P_{2n}H)(z) = \sum_{j=0}^{2n+m}{a_{j,n} z^j} \qquad \text{\rm and} \qquad (P_{2n}G)(z) = \sum_{j=0}^{2n+m}{b_{j,n} z^j}\,. \tag 3.3$$
(Note that $b_{n+m,n} = 0$. Due to the minimality of $H$ we have
$$\sup_{n \in {\Bbb N}}{\text {\rm NC}(P_{2n}G)} = \infty\,. \tag 3.4$$
Observe that (3.2) implies
$$a_{j,n} = b_{j-1,n} - \alpha b_{j,n}\,, \qquad j=1,2,\ldots,2n+m\,. \tag 3.5$$
Combining (3.1), (3.3), and (3.5), we can deduce that
$$\mu := \sup_{n \in {\Bbb N}}{|j:1 \leq j \leq 2n+m, b_{j-1,n} \neq \alpha b_{j,n}|} < \infty\,.$$
Let
$$A_n := \{j: 1 \leq j \leq 2n+m, b_{j-1,n} \neq \alpha b_{j,n}\} = \{j_{1,n} < j_{2,n} < \cdots < j_{u_n,n}\}\,,$$
where $u_n \leq \mu$ for each $n \in {\Bbb N}$, and let $j_{0,n}:=1$ and $j_{u_n+1,n}:=2n+m$.
As $\alpha \in {\Bbb C}$ with $|\alpha| = 1$ is not a root of unity, the inequality
$$j_{l+1,n} - j_{l,n} \geq |S^*|$$
for some $l=0,1,\ldots,u_n$ implies
$$b_{j,n} = 0, \qquad j = j_{l,n},\,j_{l,n} + 1,\,j_{l,n} + 2, \ldots,\,j_{l+1,n} - 1\,.$$
But then $b_{j,n} \neq 0$ is possible only for $(\mu+1)|S^*|$ values of $j = 1,2,\ldots,2n+m$, which contradicts (3.4).
This finishes the proof of the fact that each zero of $H$ is a root of unity.

Now we prove that each zero of $H$ is simple.
Without loss of generality it is sufficient to prove that $H(1) = 0$ implies that $H^{\prime}(1) \neq 0$,
the general case can easily be reduced to this. Assume to the contrary that $H(1) = 0$ and $H^{\prime}(1) \neq 0$.
Let $G_1 \in {\Cal P}_{m-1}^c$ and $G_2 \in {\Cal P}_{m-2}^c$ be defined by
$$G_1(z) := \frac{H(z)}{z-1} \qquad \text{\rm and} \qquad G_2(z) := \frac{H(z)}{(z-1)^2} = \frac{G_1(z)}{z-1}\,, \tag 3.6$$
respectively. Let
$$(P_{2n}H)(z) = \sum_{j=0}^{2n+m}{a_{j,n} z^j}\,, \tag 3.7$$
$$(P_{2n}G_1)(z) = \sum_{j=0}^{2n+m}{b_{j,n} z^j} \qquad \text{\rm and} \qquad (P_{2n}G_2)(z) = \sum_{j=0}^{2n+m}{c_{j,n} z^j}\,.$$
Due to the minimality of the degree of $H$ we have
$$\sup_{n \in {\Bbb N}}{\text {\rm NC}(P_{2n}G_1)} = \infty\,. \tag 3.8$$
Observe that (3.6) implies
$$a_{j,n} = b_{j-1,n} - b_{j,n}\,, \qquad j=1,2,\ldots,2n+m\,. \tag 3.9$$
and
$$b_{j,n} = c_{j-1,n} - c_{j,n}\, \qquad j=1,2,\ldots,2n+m\,. \tag 3.10$$
Combining (3.1), (3.7), and (3.9), we can deduce that
$$\mu := \sup_{n \in {\Bbb N}}{|j: 1 \leq j \leq 2n+m, b_{j-1,n} \neq b_{j,n}|} < \infty\,. \tag 3.11$$ 
By using (3.8) and (3.11), for every $n \in {\Bbb N}$ and  $N \in {\Bbb N}$ there is an
$L \in {\Bbb N}$ such that
$$0 \neq b := b_{L,n} = b_{L+1,n} = \cdots = b_{L+N,n}\,.$$
Combining this with (3.10), we get
$$c_{j-1,n} = c_{j,n} - b\,, \quad j = L+1,L+2,\ldots,L+N\,.$$
Hence
$$\sup_{n \in {\Bbb N}}{\max \{|c_{j,n}|: \enskip j = 0,1,\ldots,2n+m\}} = \infty\,. \tag 3.12$$
On the other hand $P_{2n} \in {\Cal P}_{2n}^c(S)$ together with the fact that the set
$S$ is finite implies that the set
$$\{|c_{j,n}|: \enskip j = 0,1,\ldots,2n+m, \enskip n \in {\Bbb N}\}$$
is also finite. This contradicts (3.12), and the proof of the fact that each zero of $H$ is simple
is finished.
\qed \enddemo

\proclaim{Lemma 3.2}
If $S \subset {\Bbb C}$ is a finite set, $P_{2n} \in {\Cal P}_{2n}^c(S)$, $H(z) := z^k - 1$, 
$$\mu := \sup_{n \in {\Bbb N}}{\text {\rm NC}(P_{2n}H)} < \infty\,, \tag 3.13$$
then there are constants $c_1 > 0$ and $c_2 > 0$ depending only on $\mu$, $k$, and $S$ 
and independent of $n$ and $\delta$ such that
$$\int_{-\delta}^{\delta}{|P_{2n}(e^{it})| \, dt} > c_1 \log \text {\rm NC}_k(P_{2n} - c_2 \delta^{-1}$$
for every $\delta \in (0,\pi)$, and hence assumption (2.1) implies 
$$\lim_{n \rightarrow \infty}{\int_{-\delta}^{\delta}|P_{2n}(e^{it})| \, dt} = \infty$$
for every $\delta \in (0,\pi)$.
\endproclaim

\demo{Proof}
We define
$$G(z) := \sum_{j=0}^{k-1}{z^j}$$
so that $H(z) = G(z)(z-1)$. Let $S_n^*$ be the set of the coefficients of $P_{2n}G$. We define
$$S^* := \bigcup_{n=1}^\infty{S_n^*}\,.$$ 
As $P_{2n} \in {\Cal P}_{2n}^c(S)$ and the set $S$ is finite, the set $S^*$ is also finite.
So by Theorem 1.3 there is an absolute constant $c > 0$ such that
$$\int_0^{2\pi}{|(P_{2n}G)(e^{it})| \, dt} \geq c\gamma \log(\text {\rm NC}(P_{2n}G)) 
\geq c\gamma \log(\text {\rm NC}_k(P_{2n}) \qquad n \in {\Bbb N} \tag 3.14$$
with
$$\gamma := \min_{z \in S^* \setminus \{0\}}{|z|}\,.$$
Observe that
$$|(P_{2n}G)(e^{it})| = \frac{1}{|e^{it}-1|} \, |(P_{2n}H)(e^{it})| \leq \frac{\mu}{|e^{it}-1|} = 
\frac{\mu}{2\sin(t/2)} \leq \frac{\pi\mu}{2t}\,, \quad t \in (-\pi,\pi)\,.$$
Hence
$$\int_{\delta}^{2\pi - \delta}{|(P_{2n}G)(e^{it})| \, dt} = 
\int_{-\pi + \delta}^{\pi - \delta}{|(P_{2n}G)(e^{it})| \, dt}\leq 2\pi \frac{\pi\mu}{2\delta} = 
\frac{\pi^2\mu}{\delta}\tag 3.15$$
Now (3.14) and (3.15) give
$$\split \int_{-\delta}^{\delta}{|P_{2n}(e^{it})| \, dt} & \geq 
\frac 1k \int_{-\delta}^{\delta}{|(P_{2n}G)(e^{it})| \, dt} \cr 
& = \frac 1k \left( \int_{0}^{2\pi}{|(P_{2n}G)(e^{it})| \, dt} - 
\int_{\delta}^{2\pi - \delta}{|(P_{2n}G)(e^{it})| \, dt} \right) \cr  
& \geq \frac 1k \, c\gamma \log(\text {\rm NC}_k(P_{2n})  - \frac{\pi^2\mu}{\delta}\,. \cr \endsplit$$
\qed \enddemo

\proclaim{Lemma 3.3}
If $S \subset {\Bbb R}$ is a finite set, $P_{2n} \in {\Cal P}_{2n}^c(S)$ are self-reciprocal, 
$H(z) := z^k - 1$, (3.13) holds, 
$$T_n(t) := P_{2n}(e^{it})e^{-int} \qquad and \qquad R_n(x) := \int_0^x{T_n(t) \, dt}\,,$$
and $0 < \delta \leq (2k)^{-1}$, then
$$\sup_{n \in {\Bbb N}} {\max_{x \in [-\delta,\delta]}{|R_n(x)|}} < \infty\,.$$
\endproclaim

\demo{Proof}
Let
$$T_n(t) = a_{0,n} + \sum_{j=1}^n{2a_{j,n}\cos(jt)}\,, \qquad a_{j,n} \in S\,.$$
Observe that (3.13) implies that
$$\sup_{n \in {\Bbb N}}{|\{j: k \leq j \leq n, a_{j-k,n} \neq a_{j,n}\}|} < \infty\,. \tag 3.16$$ 
We have
$$R_n(x) = a_{0,n}x + \sum_{j=1}^n{\frac{2a_{j,n}}{j} \, \sin(jx)}\,.$$
Now (3.16) implies that
$$R_n(x) = a_{0,n}x + \sum_{m=1}^{u_n}{F_{m,n}(x)}\,,$$
where
$$F_{m,n}(x) := \sum_{j=0}^{n_m}{\frac{2A_m \sin((j_m+jk)x)}{j_m + jk}}$$
with some $A_m \in S, m=1,2,\ldots,u_n$, and $j_m \in {\Bbb N}$, where
$$\mu := \sup_{n \in {\Bbb N}}{u_n} < \infty\,.$$
Hence it is sufficient to prove that
$$\sup_{n \in {\Bbb N}}{\max_{x \in [-\delta,\delta]}{|F_{m,n}(x)|}} < \infty\,, \qquad m=1,2,\ldots,u_n\,. \tag 3.17$$
Let $x \in (0,\delta]$, where $0 < \delta \leq (2k)^{-1}$. We break the sum as
$$F_{m,n} = 2A_m(R_{m,n} + S_{m,n})\,, \tag 3.18$$
where
$$R_{m,n}(x) := \sum_{j=0 \atop j_m + jk \leq x^{-1}}^{n_m}{\frac{\sin((j_m+jk)x)}{j_m + jk}}$$
and
$$S_{m,n}(x) := \sum_{j=0 \atop x^{-1} < j_m + jk}^{n_m}{\frac{\sin((j_m+jk)x)}{j_m + jk}}\,.$$
Here
$$|R_{m,n}(x)| \leq \sum_{j=0 \atop j_m + jk \leq x^{-1}}{\left| \frac{\sin((j_m+jk)x)}{j_m + jk} \right|} \leq |x^{-1}||x| \leq 1\,. \tag 3.19$$
Further, using Abel rearrangement, we have
$$\split S_{m,n}(x) = & -\frac{B_{m,v,k}(x)}{j_m + vk} + \frac{B_{m,u,k}(x)}{j_m + n_mk} \cr
& + \sum_{j=0 \atop x^{-1} \leq j_m + jk}^{n_m}{B_{m,j,k}(x)}\left(\frac{1}{j_m+jk} - \frac{1}{j_m+(j+1)k} \right)\,. \cr \endsplit$$ 
with
$$B_{m,j,k}(x) := \sum_{\alpha = 0}^j{\sin((j_m + \alpha k)x)}$$
and with some $u,v \in {\Bbb N}_0$ for which $x^{-1} < j_m + (u+1)k$ and $x^{-1} < j_m + (v+1)k$. 
Hence,
$$\split |S_{m,n}(x)| \leq & \left( \left| \frac{B_{m,v,k}(x)}{j_m + vk} \right| + \left| \frac{B_{m,u,k}(x)}{j_m + uk} \right| \right) \cr 
& + \sum_{j=0 \atop x^{-1} < j_m + jk}^{n_m}{|B_{m,j,k}(x)}|\left(\frac{1}{j_m+jk} - \frac{1}{j_m+(j+1)k} \right)\,. \cr \endsplit \tag 3.20$$
Observe that 
$$x^{-1} < j_m + (w+1)k < 2(j_m + wk) \qquad \text {\rm if} \enskip w \geq 1\,,$$
and
$$2k \leq \delta^{-1} \leq x^{-1} < j_m + k \qquad \text {\rm if} \enskip w = 0\,,$$
hence, 
$$\frac{1}{j_m + wk} \leq 2x\,, \qquad w \in {\Bbb N}_0 \tag 3.21$$
Observe that $x \in (0,\delta]$ and $0 < \delta \leq (2k)^{-1}$ imply that $0 < x < \pi k^{-1}$. Hence, with $z = e^{ix}$ we have
$$\split |B_{m,j,k}(x)| & = \left|\frac 12 \text {\rm Im} \left(\sum_{\alpha = 0}^j{z^{j_m + \alpha k}} \right) \right|
\leq  \left|\frac 12 \sum_{\alpha = 0}^j{z^{j_m + \alpha k}}\right| = \left|\frac 12 \sum_{\alpha = 0}^j{z^{\alpha k}}\right| \cr 
& = \left| \frac 12 \, \frac{1-z^{(j+1)k}}{1-z^k} \right| \leq \frac 12 \, |1-z^{(j+1)k}| \, \frac{1}{|1-z^k|} \leq \frac{1}{|1-z^k|} \cr 
& \leq \frac{1}{\sin(kx/2)} \leq \frac{\pi}{kx}\,. \cr \endsplit \tag 3.22$$
Combining (3.20), (3.21), and (3.22), we conclude
$$|S_{m,n}(x)| \leq \frac{\pi}{kx} \, 2x +  \frac{\pi}{kx} \, 2x +  \frac{\pi}{kx} \, 2x \leq \frac{6\pi}{k}\,. \tag 3.23$$
As $F_{n,m}$ is odd, (3.18), (3.19), and (3.23) give (3.17).
\qed \enddemo

Our next lemma is well known and may be proved simply by contradiction.

\proclaim {Lemma 3.4}
If $R$ is a continuously differentiable function on the interval $[-\delta,\delta]$, $\delta > 0$,
$$\int_{-\delta}^{\delta}{|R^{\prime}(x)| \, dx} = L \qquad \text {\rm and} \qquad \max_{x \in [-\delta,\delta]}{|R(x)|} = M\,,$$
then there is an $\eta \in [-M,M]$ such that $R-\eta$ has at least $L(2M)^{-1}$ zeros in $[-\delta,\delta]$.
\endproclaim

\proclaim {Lemma 3.5}
If $S \subset {\Bbb R}$ is a finite set, $P_{2n} \in {\Cal P}_{2n}^c(S)$ are self-reciprocal, 
$$T_n(t) := P_{2n}(e^{it})e^{-int}\,,$$
$H(z) := z^k - 1$, and  (3.13) and (2.1) hold, then (2.2) also holds. 
\endproclaim

\demo{Proof}
Let $0 < \delta \leq (2k)^{-1}$. Let $R_n$ be defined by
$$R_n(x) := \int_0^x{T_n(t) \, dt}\,.$$
Observe that $|T_n(x)| = |P_{2n}(e^{ix})|$ for all $x \in {\Bbb R}$.
By Lemmas 3.2 and 3.3 we have
$$\lim_{n \rightarrow \infty} {\int_{-\delta}^{\delta}{|R_n^\prime(x)| \, dx}} =
\lim_{n \rightarrow \infty} {\int_{-\delta}^{\delta}{|T_n(x)| \, dx}} =
\lim_{n \rightarrow \infty} {\int_{-\delta}^{\delta}{|P_{2n}(e^{ix})| \, dx}} =
\infty$$
and
$$\sup_{n \in {\Bbb N}}{\max_{[-\delta,\delta]}{|R_n(x)|}} < \infty\,.$$
Therefore, by Lemma 3.4 there are $c_n \in {\Bbb R}$ such that
$$\lim_{n \rightarrow \infty}{\text {\rm NZ}(R_n - c_n)} = \infty\,.$$
However, $T_n(x) = (R_n - c_n)^{\prime}(x)$ for all $x \in {\Bbb R}$, and hence
$$\lim_{n \rightarrow \infty}{\text {\rm NZ}(T_n)} = \infty\,.$$
\qed \enddemo

Our next lemma follows immediately from Lemmas 3.1 and 3.5. 

\proclaim {Lemma 3.6}
If $S \subset {\Bbb R}$ is a finite set, $P_{2n} \in {\Cal P}_{2n}^c(S)$ are self-reciprocal, 
$$T_n(t) := P_{2n}(e^{it})e^{-int}\,,$$ 
(2.1) holds, and there is a polynomial $H \in {\Cal P}_m$ 
such that (3.13) holds, then (2.2) also holds.  
\endproclaim

Moreover, we have the following two observation.

\proclaim {Lemma 3.7}
Let $(n_\nu)$ be a strictly increasing sequence of positive integers.
If $S \subset {\Bbb R}$ is a finite set, $P_{2n_\nu} \in {\Cal P}_{2n_\nu}^c(S)$ 
are self-reciprocal, $T_{n_\nu}(t) := P_{2n_\nu}(e^{it})e^{-in_\nu t}$,
$$\lim_{\mu \rightarrow \infty}{\text {\rm NC}_k(P_{2n_\mu})} = \infty$$
for every $k \in {\Bbb N}$, and there is a polynomial $H \in {\Cal P}_m$ such that 
$$\sup_{\nu \in {\Bbb N}}{\text {\rm NC}(P_{2n_\nu}H)} < \infty\,,$$
then
$$\lim_{\nu \rightarrow \infty}{\text {\rm NZ}(T_{n_\nu})} = \infty\,.$$
\endproclaim 

\demo{Proof}
We define 
$$P_{2n} := P_{2n_\nu}\,, \qquad n_\nu \leq n < n_{\nu+1}\,,$$
and apply Lemma 3.6.  
\qed \enddemo

The next lemma is straightforward consequences of Theorem 1.4. 

\proclaim{Lemma 3.8} Let $\lambda_0 < \lambda_1 < \cdots < \lambda_m$
be nonnegative integers and let
$$Q_m(t) = \sum_{j=0}^m{A_j \cos(\lambda_jt)}\,, \qquad A_j \in {\Bbb R}\,,
\enskip j = 0,1,\ldots, m\,.$$
Then
$$\int_{-\pi}^{\pi}{|Q_m(t)|\,dt} \geq \frac{1}{60}
\sum_{j=0}^m{\frac{|A_{m-j}|}{j+1}}\,.$$
\endproclaim

We will also need the lemma below in the proof of Theorem 2.1..

\proclaim{Lemma 3.9} Let $\lambda_0 < \lambda_1 < \cdots < \lambda_m$
be nonnegative integers and let
$$Q_m(t) = \sum_{j=0}^m{A_j \cos(\lambda_jt)}\,, \qquad A_j \in {\Bbb R}\,,
\enskip j = 0,1,\ldots, m\,.$$
Let $A := \max\{|A_j|:j=0,1,\ldots, m\}\,.$
Suppose $Q_m$ has at most $K-1$ zeros in the period $[-\pi,\pi)$.
Then
$$\int_{-\pi}^{\pi}{|Q_m(t)|\,dt} \leq 2KA \left(\pi + \sum_{j=1}^m{\frac{1}{\lambda_j}}\right)
\leq 2KA(5 + \log m)\,.$$
\endproclaim

\demo{Proof} We may assume that $\lambda_0 = 0$, the case $\lambda_0 > 0$ can be handled
similarly. Associated with $Q_m$ in the lemma let
$$R_m(t) := A_0t + \sum_{j=0}^m{\frac{A_j}{\lambda_j} \sin(\lambda_jt)}\,.$$
Clearly
$$\max_{t \in [-\pi,\pi]}{|R_m(t)|} \leq A\left( \pi + \sum_{j=1}^m{\frac{1}{\lambda_j}} \right)\,.$$
Also, for every $c \in {\Bbb R}$ the function $R_m - c$ has at
most $K$ zeros in the period $[-\pi,\pi)$, otherwise Rolle's Theorem implies
that $Q_m = (R_m - c)^\prime$ has at least $K$ zeros in the period $[-\pi,\pi)$.
Hence
$$\split \int_{-\pi}^{\pi}{|Q_m(t)|\,dt} & = \int_{-\pi}^{\pi}{|R_m^{\prime}(t)|\,dt} = 
V_{-\pi}^\pi(R_m) \leq 2K \max_{t \in [-\pi,\pi]}{|R_m(t)|} \cr
& \leq 2KA \left(\pi + \sum_{j=1}^m{\frac{1}{\lambda_j}}\right) 
\leq 2KA(5 + \log m)\,, \cr \endsplit$$
and the lemma is proved. 
\qed \enddemo 

\proclaim{Lemma 3.10} Suppose $k \in {\Bbb N}$.
Let
$$z_j := \exp\left(\frac{2\pi ji}{k}\right)\,, \qquad j=0,1,\ldots, k-1\,,$$
be the $k$th roots of unity. Suppose
$$0 \notin \{b_0,b_1, \ldots, b_{k-1}\} \subset {\Bbb R}$$
and
$$Q(z) := \sum_{j=0}^{k-1}{b_jz^j}\,.$$
Then there is a value of $j \in \{0,1,\ldots\, k-1\}$
for which $\text{\rm Re}(Q(z_j)) \neq 0\,$.
\endproclaim

\demo{Proof} If the statement of the lemma were false, then
$$z^{k-1}(Q(z)+Q(1/z))=(z^k-1) \sum_{\nu=0}^{k-2}{\alpha_{\nu}z^{\nu}}\,.$$
Observe that the coefficient of $z^{k-1}$ on the right hand side is $0$, 
while the coefficient of $z^{k-1}$ on the left hand side is $b_1 \neq 0$, 
a contradiction. 
\qed \enddemo

\proclaim{Lemma 3.11} If $0 \notin \{b_0,b_1, \ldots, b_{k-1}\} \subset {\Bbb R}$,
$\{a_0, a_1, \ldots, a_{m-1}\} \subset {\Bbb R}$, where $m = uk$ with some integer $u \geq 0$,
$$a_{m+lk+j} = b_j\,, \qquad l=0,1,\ldots\,, \quad j=0,1,\ldots, k-1\,,$$
and $n=m+lk+r$ with integers $m \geq 0$, $l\geq 0$, $k\geq 1$, and $0\leq r \leq k-1$, then 
there is a constant $c_9 > 0$ independent of $n$ such that
$$T_n(t) := \text{\rm Re}\left(\sum_{j=0}^n{a_je^{ijt}}\right)$$
has at least $c_9n$ zeros in $[-\pi,\pi)$.
\endproclaim

\demo{Proof} Note that
$$\sum_{j=0}^n{a_jz^j} = \sum_{j=0}^{m-1}{a_jz^j} +
z^m\left(\sum_{j=0}^{k-1}{b_jz^j}\right)\frac{z^{(l+1)k}-1}{z^k-1}
+z^{m+lk}\sum_{j=0}^r{b_jz^j} = P_1(z) + P_2(z)\,,$$
where
$$P_1(z) := \sum_{j=0}^{m-1}{a_jz^j} + z^{m+lk}\sum_{j=0}^u{b_jz^j}$$
and
$$P_2(z) := z^{uk} \sum_{j=0}^{k-1}{b_jz^j} \frac{z^{(l+1)k}-1}{z^k-1}
= Q(z) z^{uk}\frac{z^{(l+1)k}-1}{z^k-1}$$
with
$$Q(z) := \sum_{j=0}^{k-1}{b_jz^j}\,.$$
By Lemma 3.5 there is a $k$th root of unity $\xi = e^{i\tau}$
such that $\text{Re}(Q(\xi)) \neq 0$. Then, for every $K > 0$
there is a $\delta \in (0,2\pi/k)$ such that $\text{Re}(P_2(e^{it}))$
oscillates between $-K$ and $K$ at least $c_{10}(l+1)k\delta$ times,
where $c_{10} > 0$ is a constant independent of $n$. Now we choose $\delta \in (0,2\pi/k)$
for
$$K := 1 + \sum_{j=0}^{m-1}{|a_j|} + \sum_{j=0}^{k-1}{|b_j|}\,.$$
Then
$$T_n(t) := \text{Re}\left(\sum_{j=0}^n{a_je^{ijt}}\right)
= \text{Re}(P_1(e^{it})) + \text{Re}(P_2(e^{it}))$$
has at least one zero on each interval on which $\text{Re}(P_2(e^{it}))$ oscillates between $-K$ and $K$, and
hence it has at least $c_{10}(l+1)k\delta > c_9n$ zeros on $[-\pi,\pi)$,
where $c_9 > 0$ is a constant independent of $n$.
\qed \enddemo

\head Proof of the Theorems \endhead

We denote the set of all real trigonometric polynomials of degree at most $k$ by 
${\Cal T}_k$.

\demo{Proof of Theorems 2.1}
Suppose the theorem is false. Then there are $k \in {\Bbb N}$, a strictly increasing sequence
$(n_\nu)_{\nu = 1}^\infty$ of positive integers and
even trigonometric polynomials $Q_{n_\nu} \in {\Cal T}_k$ with maximum norm $1$
on the period such that
$$T_{n_\nu}(t)Q_{n_\nu}(t) \geq 0\,, \qquad t \in {\Bbb R}\,. \tag 4.1$$
We can pick a subsequence of $(n_\nu)_{\nu = 1}^\infty$ (without loss of
generality we may assume that it is the sequence $(n_\nu)_{\nu = 1}^\infty$ itself)
that converges to a $Q \in {\Cal T}_k$ uniformly on the period $[-\pi,\pi)$.
That is,
$$\lim_{\nu \rightarrow \infty}{\varepsilon_\nu} = 0 \qquad \text{with} \qquad 
\varepsilon_\nu := \max_{t \in [-\pi,\pi]}{|Q(t) - Q_{n_\nu}(t)|}\,. \tag 4.2$$
We introduce the notation 
$$T_{n_\nu}(t)Q(t)^3 = \left( \sum_{j=0}^{n_\nu}{a_{j,\nu} \cos(jt)} \right) Q(t)^3  
= \sum_{j=0}^{K_\nu}{b_{j,\nu}\cos(\beta_{j,\nu}t)}\,, \tag 4.3$$
$$b_{j\nu} \neq 0, \qquad j = 0,1,\ldots,K_{\nu}\,,$$
and
$$T_{n_\nu}(t)Q(t)^4 = \left( \sum_{j=0}^{n_\nu}{a_{j,\nu} \cos(jt)} \right) Q(t)^4 
= \sum_{j=0}^{L_\nu}{d_{j,\nu} \cos(\delta_{j,\nu}t)}\,, \tag 4.4$$
$$d_{j,\nu} \neq 0\,, \qquad j=0,1,\ldots,L_{\nu}\,,$$
where $\beta_{0,\nu} < \beta_{1,\nu} < \cdots < \beta_{K_\nu,\nu}$ and 
$\delta_{0,\nu} < \delta_{1,\nu} < \cdots < \delta_{L_\nu,\nu}$ 
are nonnegative integers. Since the set 
$\{a_{j,\nu}: j=0,1,\ldots,n_\nu, \enskip \nu \in {\Bbb N}\} \subset {\Bbb R}$ is finite, the sets
$$\{b_{j,\nu}: j=0,1,\ldots,K_\nu, \enskip \nu \in {\Bbb N}\} \subset {\Bbb R} \qquad \text{and} \qquad 
\{d_{j,\nu}: j=0,1,\ldots,L_\nu, \enskip \nu \in {\Bbb N}\} \subset {\Bbb R}$$
are finite as well. Hence there are $\rho, M \in (0,\infty)$ such that
$$|a_{j,\nu}| \leq M\,, \qquad j=0,1,\ldots,n_\nu, \enskip \nu \in {\Bbb N}\,, \tag 4.5$$
$$\rho \leq |b_{j,\nu}|\,, \qquad j=0,1,\ldots,K_\nu, \enskip \nu \in {\Bbb N}\,, \tag 4.6$$
and     
$$|d_{j,\nu}| \leq M\,, \qquad j=0,1,\ldots,L_\nu, \enskip \nu \in {\Bbb N}\,. \tag 4.7$$
Observe that our indirect assumption together with Lemma 3.7 implies that
$$\lim_{\nu \rightarrow \infty}{K_\nu} = \infty \qquad \text{and} \qquad 
\lim_{\nu \rightarrow \infty}{L_\nu} = \infty\,. \tag 4.8$$
We claim that
$$K_\nu \leq c_3L_\nu \tag 4.9$$
with some $c_3 > 0$ independent of $\nu \in {\Bbb N}$. Indeed, using Parseval's formula (4.2), (4.3), and (4.6)
we deduce
$$\frac{1}{\pi} \, \int_{-\pi}^\pi {T_{n_\nu}(t)^2Q(t)^4Q_{n_\nu}(t)^2\,dt} =
\frac{1}{\pi} \, \int_{-\pi}^\pi{(T_{n_\nu}(t)Q(t)^2Q_{n_\nu}(t))^2\,dt}  \geq \frac 12 \rho^2 K_\nu \tag 4.10$$
for every sufficiently large $\nu \in {\Bbb N}\,.$
Also, (4.1) -- (4.8) imply
$$\split & \frac{1}{\pi} \, \int_{-\pi}^\pi {T_{n_\nu}(t)^2 Q(t)^4 Q_{n_\nu}(t)^2 \,dt} =
\frac{1}{\pi} \, \int_{-\pi}^\pi {(T_{n_\nu}(t)Q_{n_\nu}(t))(T_{n_\nu}(t)Q(t)^4) Q_{n_\nu}(t) \,dt} \cr 
\leq & \, \frac{1}{\pi} \, \left( \int_{-\pi}^\pi {T_{n_\nu}(t)Q_{n_\nu}(t) \,dt} \right)
\left( \max_{t \in [-\pi,\pi]} {|T_{n_\nu}(t)Q(t)^4|} \right) 
\left( \max_{t \in [-\pi,\pi]}{|Q_{n_\nu}(t)|} \right) \cr
\leq & \,  \frac{1}{\pi} \, \left( \int_{-\pi}^\pi {T_{n_\nu}(t)Q_{n_\nu}(t) \,dt} \right) 
L_\nu M \left( \max_{t \in [-\pi,\pi]}{|Q_{n_\nu}(t)|} \right) \cr 
\leq & \, c_4  L_\nu \cr \endsplit \tag 4.11$$
with a constant $c_4 > 0$ independent of $\nu$ for every $\nu \in {\Bbb N}$.
Now (4.9) follows from (4.10) and (4.11). From Lemma 3.8 we deduce
$$\int_{-\pi}^\pi {|T_{n_\nu}(t)Q(t)^4| \,dt} \geq c_5 \rho \log L_\nu \tag 4.12$$
with some constant $c_5 > 0$ independent of $\nu \in {\Bbb N}\,.$
On the other hand, using (4.1), Lemma 3.9, (4.2), (4.4), (4.9), and (4.8), we obtain
$$\split & \int_{-\pi}^\pi {|T_{n_\nu}(t) Q(t)^4|\,dt} \cr
\leq & \int_{-\pi}^\pi {|T_{n_\nu}(t)Q(t)^3||Q_{n_\nu}(t)|\,dt}
+ \int_{-\pi}^\pi {|T_{n_\nu}(t)Q(t)^3|\,|Q(t) - Q_{n_\nu}(t)|\,dt} \cr 
\leq & \int_{-\pi}^\pi {|T_{n_\nu}(t)Q_{n_\nu}(t)||Q(t)^3|\,dt}
+ \int_{-\pi}^\pi {|T_{n_\nu}(t)Q(t)^3|\,|Q(t) - Q_{n_\nu}(t)|\,dt} \cr
\leq & \left( \int_{-\pi}^\pi {|T_{n_\nu}(t)Q_{n_\nu}(t)|\,dt} \right)
\left( \max_{t \in [-\pi,\pi]}{|Q(t)|^3} \right) \cr
& \qquad \qquad \qquad + \left( \int_{-\pi}^\pi {|T_{n_\nu}(t)Q(t)^3|\,dt} \right) 
\left( \max_{t \in [-\pi,\pi]}{|Q(t) - Q_{n_\nu}(t)|} \right) \cr
\leq & \left( \int_{-\pi}^\pi {T_{n_\nu}(t)Q_{n_\nu}(t)\,dt} \right)
\left( \max_{t \in [-\pi,\pi]}{|Q(t)|^3} \right) + \left( \int_{-\pi}^\pi {|T_{n_\nu}(t)Q(t)^3|\,dt} \right)
\varepsilon_\nu \cr
\leq & c_6 + c_7(\log K_\nu) \varepsilon_\nu
\leq c_6 + c_7(\log (c_3 L_\nu)) \varepsilon_\nu \cr
\leq & c_8 + c_7(\log L_\nu) \varepsilon_\nu = o( \log L_\nu )\,, \cr 
\endsplit \tag 4.13$$
where $c_6, c_7$, and $c_8$ are constants independent of $\nu \in {\Bbb N}$.
Since (4.13) contradicts (4.12), the proof of the theorem is finished.
\qed \enddemo

\demo{Proof of Corollary 2.2}
Observe that assumption (2.3) implies assumption (2.1). 
\qed \enddemo

\demo{Proof of Corollary 2.3}
Corollary 2.2 implies
$$\lim_{k \rightarrow \infty}{\text {\rm NZ}(P_{2k})} = \infty\,. \tag 4.14$$
and 
$$\lim_{k \rightarrow \infty}{\text {\rm NZ}(P_{2k+1})} = \infty\,. \tag 4.15$$
Note that (4.14) is an obvious consequence of Theorem 2.1. To see (4.15) 
observe that if $P_{2k+1} \in {\Cal P}_{2k+1}^c(S)$ are self-reciprocal then 
$P_{2k+2}$ defined by 
$$\widetilde{P}_{2k+2}(z) := (z+1)P_{2k+1}(z)\in {\Cal P}_{2k+2}^c(\widetilde{S})$$ 
are also self-reciprocal, where the finiteness of $S$ implies the finiteness of 
$\widetilde{S}$. Also 
$$\lim_{n \rightarrow \infty}{|P_n(1)|} = \infty\,.$$ 
implies 
$$\lim_{k \rightarrow \infty}{|\widetilde{P}_{2k+2}(1)|} = 
\lim_{k \rightarrow \infty}{2|P_{2k+1}(1)|} = \infty\,.$$
Hence the polynomials $\widetilde{P}_{2k+2} \in {\Cal P}_{2k+2}^c(\widetilde{S})$ 
satisfy the assumptions of Corollary 2.2. 
\qed \enddemo

\demo{Proof of Corollary 2.4}
If the finite set $S \subset {\Bbb R}$ has property (2.6), then assumption (2.1) 
is satisfied.  
\qed \enddemo

\demo{Proof of Corollary 2.5}
Corollary 2.4 implies (4.14) and (4.15).
Note that (4.14) is an obvious consequence of Corollary 2.4. To see (4.15)
observe that if $P_{2k+1} \in {\Cal P}_{2k+1}^c(S)$ are self-reciprocal then $P_{2k+2}$ 
defined by
$$\widetilde{P}_{2k+2}(z) := (z+1)P_{2k+1}(z)\in {\Cal P}_{2k+2}^c(\widetilde{S})$$
are also self-reciprocal, where the finiteness of $S$ implies the finiteness of
$$\widetilde{S}:=\{s_1 + s_2: \enskip s_1, s_2 \in S \cup \{0\}\}\,.$$ 
Also, it is easy to see that if
$$s_1 + s_2 + \cdots + s_k = 0$$
with some $s_1,s_2,\ldots,s_k \in \widetilde{S}$ and $k \in {\Bbb N}$, then 
$s_1=s_2= \cdots =s_k=0$.
Hence the polynomials $\widetilde{P}_{2k+2} \in {\Cal P}_{2k+2}^c(\widetilde{S})$
satisfy the assumptions of Corollary 2.4.
\qed \enddemo

\demo{Proof of Corollary 2.6}
This is an obvious consequence of Corollary 2.2.
\qed \enddemo

\demo{Proof of Corollary 2.6}
This is an obvious consequence of Corollary 2.6.
\qed \enddemo

\head 14. Acknowledgements \endhead
The author wishes to thank Stephen Choi and Jonas Jankauskas for their reading earlier versions of 
my paper carefully, pointing out many misprints, and their suggestions to make the paper more readable.

\enddocument